\documentclass[11pt]{amsart}
\usepackage{amssymb,amsmath,amsthm}
\usepackage{booktabs}
\usepackage[margin=1.15in]{geometry}

\newtheorem{thm}{Theorem}[section]
\newtheorem{lem}[thm]{Lemma}
\newtheorem{prop}[thm]{Proposition}
\newtheorem{cor}[thm]{Corollary}

\theoremstyle{definition}
\newtheorem{defn}[thm]{Definition}
\newtheorem{example}[thm]{Example}
\newtheorem{rem}[thm]{Remark}
\newtheorem*{dietz}{Dietzmann's Lemma}

\newtheorem*{thm33}{Theorem 3.3}
\newtheorem*{cor34}{Corollary 3.4}
\newtheorem*{cor35}{Corollary 3.5}
\newtheorem*{cor41}{Corollary 4.1}

\newcommand{\bt}{\circledast}   % general tensor product (Definition 1.5)
\newcommand{\boxt}{\boxtimes}  % box-tensor product of [LT12]
\DeclareMathOperator{\Aut}{Aut}
\DeclareMathOperator{\Sym}{Sym}
\DeclareMathOperator{\coker}{coker}

\title{Tensor products of finite groups without compatibility conditions are finite}
\author{Sathasivam Kalithasan, Viji Zachariah Thomas}
\subjclass[2010]{20F05, 20F24, 20F99, 20J99}
\keywords{Nonabelian tensor product, Inassaridze's tensor product, box-tensor product, general tensor product, nonabelian homology groups, finite groups}

\begin{document}

\begin{abstract}
It has been an open question whether Inassaridze's tensor product $G\otimes H$ is finite whenever $G$ and $H$ are finite. In this paper, we answer this question in the affirmative, and, as a consequence, we prove that the low dimensional nonabelian homology groups $H_i(G,A)$, $i=0,1$, are finite whenever $G$ and $A$ are finite. We achieve this by introducing a much more general tensor product $G\bt H$, defined for groups acting on themselves and on each other by arbitrary automorphisms, with no compatibility conditions, and subject only to the two defining relations of the nonabelian tensor product. We prove that $G\bt H$ is finite whenever $G$ and $H$ are finite. Since the nonabelian tensor product of Brown and Loday, Inassaridze's tensor product and the box-tensor product are all quotients or special cases of $G\bt H$, this brings the finiteness results for all these tensor products under one umbrella.
\end{abstract}

\maketitle

\section{Introduction}

R. Brown and J.-L. Loday introduced the nonabelian tensor product $G\otimes H$ for a pair of groups $G$ and $H$ in \cite{BL84} and \cite{BL87} in the context of an application in homotopy theory. Their construction requires the mutual actions of $G$ and $H$ to satisfy certain compatibility conditions. In \cite{Ina96a}, N. Inassaridze extends this construction to arbitrary pairs of groups $G$ and $H$ which act on themselves by conjugation and each of which acts on the other, without requiring the compatibility conditions. He constructs the homology groups of groups with coefficients in any group as the left derived functors of the tensor product defined by him, thereby generalizing the classical theory of homology of groups. He also gives an application to algebraic $K$-theory of noncommutative local rings. In \cite{LT12}, Ladra and the second author introduce the box-tensor product $G\boxt H$, in which the groups are allowed to act on themselves by arbitrary automorphisms.

Consider two groups $G$ and $H$ acting on themselves and on each other. In this paper, we only consider actions by automorphisms. Let $\Aut(G)$ and $\Aut(H)$ be the automorphism groups of $G$ and $H$, respectively. Consider $\rho_G\colon G\to\Aut(G)$, $\rho_H\colon H\to\Aut(H)$, $\sigma_G\colon G\to\Aut(H)$ and $\sigma_H\colon H\to\Aut(G)$. We will write ${}^{g}h$ for $\sigma_G(g)(h)$ and ${}^{g}g'$ for $\rho_G(g)(g')$, and likewise ${}^{h}g$ for $\sigma_H(h)(g)$ and ${}^{h}h'$ for $\rho_H(h)(h')$, where $g,g'\in G$ and $h,h'\in H$. When $G$ and $H$ act on themselves by conjugation, that is ${}^{g}g'=gg'g^{-1}$ and ${}^{h}h'=hh'h^{-1}$, the free product $G\ast H$ acts on $G$ and on $H$, and for $g\in G$ and $h\in H$ we write $[g,h]=ghg^{-1}h^{-1}\in G\ast H$.

\begin{defn}[See {\cite{BL87}}]\label{def:comp}
Let $G$ and $H$ be groups that act on themselves by conjugation and each of which acts on the other. The mutual actions are said to be compatible if
\begin{align}
{}^{({}^{h}g)}h' &= {}^{hgh^{-1}}h' \qquad\text{and} \tag{1.1.1}\\
{}^{({}^{g}h)}g' &= {}^{ghg^{-1}}g', \tag{1.1.2}
\end{align}
for all $g,g'\in G$ and $h,h'\in H$.
\end{defn}

\begin{defn}[See {\cite{LT12}}]\label{def:half}
Let $G$ and $H$ be groups acting on each other and acting on themselves by conjugation. We say that the mutual actions are half compatible if either (1.1.1) or (1.1.2) holds.
\end{defn}

For groups acting on themselves by conjugation and acting on each other via arbitrary automorphisms, with no compatibility conditions, Inassaridze's tensor product is defined as follows.

\begin{defn}[See {\cite{Ina96a}}]\label{def:inas}
Let $G$ and $H$ be groups acting on themselves by conjugation and each of which acts on the other. Inassaridze's tensor product $G\otimes H$ is the group generated by the symbols $g\otimes h$ for $g\in G$ and $h\in H$ with relations
\begin{align}
gg'\otimes h &= \big({}^{g}g'\otimes{}^{g}h\big)(g\otimes h), \tag{1.3.1}\\
g\otimes hh' &= (g\otimes h)\big({}^{h}g\otimes{}^{h}h'\big), \tag{1.3.2}\\
(g\otimes h)\big(g'\otimes h'\big)(g\otimes h)^{-1} &= {}^{[g,h]}g'\otimes{}^{[g,h]}h', \tag{1.3.3}
\end{align}
for all $g,g'\in G$ and $h,h'\in H$.
\end{defn}

In \cite{Ina96a}, the author requires four relations, but as observed in \cite{LT12}, the above three relations are sufficient. If the mutual actions are compatible, then $G\otimes H$ is the nonabelian tensor product of Brown and Loday \cite{BL84,BL87}. When the mutual actions are compatible, then (1.3.1) and (1.3.2) imply (1.3.3).

\begin{defn}[See {\cite{LT12}}]\label{def:full}
Let $G$ and $H$ be groups which act on themselves and on each other. The mutual actions are said to be fully compatible if
\begin{equation}
{}^{({}^{a}b)}c={}^{a}\Big({}^{b}\big({}^{a^{-1}}c\big)\Big) \tag{1.4.1}
\end{equation}
for all $a,b,c\in G\cup H$.
\end{defn}

The following group is the main object of this paper. It is defined by the first two relations of Definition~\ref{def:inas}, for groups acting on themselves by arbitrary automorphisms.

\begin{defn}\label{def:box}
Let $G$ and $H$ be groups which act on themselves and on each other. Let $G\bt H$ be the group generated by the symbols $g\bt h$ for $g\in G$ and $h\in H$ with relations
\begin{align}
gg'\bt h &= \big({}^{g}g'\bt{}^{g}h\big)(g\bt h), \tag{1.5.1}\\
g\bt hh' &= (g\bt h)\big({}^{h}g\bt{}^{h}h'\big), \tag{1.5.2}
\end{align}
for all $g,g'\in G$ and $h,h'\in H$. We call $G\bt H$ the \emph{general tensor product} of $G$ and $H$.
\end{defn}

If the mutual actions are fully compatible, then $G\bt H$ is the box-tensor product $G\boxt H$ of \cite{LT12}. The group $G\bt H$ is also considered in \cite[Example~2.2]{LT12} for mutual actions which are not fully compatible. If $G$ and $H$ act on themselves by conjugation, then (1.5.1) and (1.5.2) are precisely (1.3.1) and (1.3.2), and hence Inassaridze's tensor product $G\otimes H$ is a quotient of $G\bt H$.

To avoid confusion, we use three different symbols throughout the paper, see Table~\ref{tab:products}. We denote the general tensor product of Definition~\ref{def:box} by $G\bt H$, and we reserve $G\boxt H$ for the box-tensor product of \cite{LT12}. Following \cite{Ina96a}, we use the same symbol $G\otimes H$ for Inassaridze's tensor product and for the nonabelian tensor product of Brown and Loday. This causes no ambiguity, since the two coincide when the mutual actions are compatible.

\begin{table}[ht]
\centering
\footnotesize
\setlength{\tabcolsep}{4pt}
\renewcommand{\arraystretch}{1.25}
\begin{tabular}{@{}llllll@{}}
\toprule
Product & Symbol & Self-actions & Mutual actions & Relations & Reference\\
\midrule
Nonabelian tensor product & $G\otimes H$ & conjugation & compatible & (1.3.1), (1.3.2) & \cite{BL84,BL87}\\
Inassaridze's tensor product & $G\otimes H$ & conjugation & arbitrary & (1.3.1)--(1.3.3) & \cite{Ina96a}\\
Box-tensor product & $G\boxt H$ & arbitrary & fully compatible & (1.5.1), (1.5.2) & \cite{LT12}\\
General tensor product & $G\bt H$ & arbitrary & arbitrary & (1.5.1), (1.5.2) & Definition~\ref{def:box}\\
\bottomrule
\end{tabular}
\medskip
\caption{The tensor products considered in this paper. All actions are by automorphisms.}
\label{tab:products}
\end{table}

In \cite{Ell87}, Ellis shows that the nonabelian tensor product of two finite groups is finite, using in his proof part of an exact sequence from homology given in \cite{BL84} and the fact that the homology of a finite group is finite. Already in \cite{BJR87} and \cite{Kap99}, the question is asked if a purely group theoretic proof can be given. Such proofs are given in \cite{Tho10} and, via the box-tensor product, in \cite{LT12}. More generally, the study of finiteness conditions for the nonabelian tensor product and its generalizations is a common theme in the literature, see for instance \cite{BNR21,BNR25} and the references therein. These results assume that the mutual actions are compatible, or fully compatible in the case of the box-tensor product. It has remained an open question whether the tensor product $G\otimes H$ introduced in \cite{Ina96a} is finite when the two groups $G$ and $H$ are finite (see \cite{LT12}). In \cite{Ina96b}, Inassaridze proves the finiteness of $G\otimes H$ when the factors are finite along with some additional hypothesis, such as the action of $H$ on $G$ being trivial together with a solvability condition \cite[Theorems~9 and~12]{Ina96b}. In \cite{LT12}, the authors generalize these results by proving that $G\otimes H$ is finite if $G$ and $H$ are finite and $G$ acts on $H$ trivially \cite[Theorem~6.3]{LT12}, and more generally if the mutual actions are half compatible \cite[Theorem~6.4]{LT12}. The purpose of this paper is to remove these additional hypotheses. In fact, we show that no compatibility condition and no relation beyond (1.5.1) and (1.5.2) is needed.

With this set up, we now state our main results.

\begin{thm33}
If $G$ and $H$ are finite groups acting on themselves and on each other, then $G\bt H$ is finite.
\end{thm33}

Since Inassaridze's tensor product is a quotient of $G\bt H$, we obtain the following.

\begin{cor34}
If $G$ and $H$ are finite groups acting on themselves by conjugation and each of which acts on the other, not necessarily compatibly, then $G\otimes H$ is finite.
\end{cor34}

Corollary~3.4 generalizes Theorems~6.3 and~6.4 of \cite{LT12}, and hence \cite[Theorems~9 and~12]{Ina96b}. Theorem~3.3 also recovers the finiteness of the box-tensor product \cite[Theorem~4.8]{LT12} and, since the nonabelian tensor product is a special case of Inassaridze's tensor product, it gives yet another proof of Ellis's theorem \cite{Ell87}.

\begin{cor35}
Let $G$ and $H$ be finite groups acting on each other fully compatibly. Then the box-tensor product $G\boxt H$ is finite. Let $G$ and $H$ be finite groups acting on themselves by conjugation and on each other compatibly. Then the nonabelian tensor product $G\otimes H$ is finite.
\end{cor35}

As a result of Corollary~3.4, we obtain the finiteness of the low dimensional nonabelian homology groups, which generalizes \cite[Corollary~6.5]{LT12}.

\begin{cor41}
Let $G$ and $A$ be finite groups acting on each other. Then the nonabelian homology groups $H_i(G,A)$ are finite when $i=0,1$.
\end{cor41}

\subsection*{Strategy of the proof}

Our main tool is Lemma~2.2, a finiteness criterion inspired by Dietzmann's Lemma and Schur's theorem: a group generated by a finite normal set is finite if and only if its abelianization is finite. Thus, to prove that $G\bt H$ is finite, it suffices to establish two facts. First, the set $Y=\{g\bt h\mid g\in G,\ h\in H\}$ of tensors is a finite normal set. Second, the abelianization $(G\bt H)^{ab}$ is finite. The paper is organized as follows: In Section~2, we collect some preliminary results. In Section~3, we prove our main theorem and deduce the finiteness of Inassaridze's tensor product, and we give an example of mutual actions which are not half compatible. In Section~4, we deduce the finiteness of the low dimensional nonabelian homology groups.

\section{Preliminaries}

We start with the following classical theorem of Schur (see \cite[10.1.4]{Rob82}).

\begin{thm}[Schur]\label{thm:schur}
If the center $Z(K)$ of a group $K$ has finite index in $K$, then the derived subgroup $[K,K]$ is finite.
\end{thm}

Noting that a subset of a group is normal if it contains all conjugates of its elements, Dietzmann's Lemma \cite{Die37} (see also \cite[14.5.7]{Rob82}) can be stated as follows.

\begin{dietz}
In any group, a finite normal subset consisting of elements of finite order generates a finite normal subgroup.
\end{dietz}

The following finiteness criterion might be well known in the study of FC-groups, but we include a proof for the convenience of the readers.

\begin{lem}\label{lem:B}
Let $K$ be a group generated by a finite normal set $Y$. Then the following are equivalent:
\begin{enumerate}
\item[(i)] $K$ is finite;
\item[(ii)] the abelianization $K^{ab}$ is finite;
\item[(iii)] the abelianization $K^{ab}$ is a torsion group.
\end{enumerate}
\end{lem}

\begin{proof}
The implications (i)$\Rightarrow$(ii)$\Rightarrow$(iii) are clear. Assume (iii). Write $Y=\{y_1,\dots,y_n\}$, and let $\bar{y}_i$ denote the image of $y_i$ in $K^{ab}=K/[K,K]$. Since $Y$ generates $K$, the elements $\bar{y}_1,\dots,\bar{y}_n$ generate $K^{ab}$. By (iii), each $\bar{y}_i$ has finite order $d_i$. Thus $|K^{ab}|\le d_1\cdots d_n$.

Since $Y$ is a normal set, $K$ acts on $Y$ by conjugation, and we obtain a homomorphism $\theta\colon K\to\Sym(Y)$. Since $Y$ generates $K$, we have $\ker\theta=\bigcap_{y\in Y}C_K(y)=Z(K)$. Thus $K/Z(K)$ embeds in $\Sym(Y)$, and so $[K:Z(K)]\le|Y|!$ is finite. By Theorem~\ref{thm:schur}, the derived subgroup $[K,K]$ is finite. Therefore $|K|=|[K,K]|\cdot|K^{ab}|$ is finite.
\end{proof}

The next lemma collects some basic properties of $G\bt H$. It is the analogue of the relation \cite[(3.4.1)]{LT12}.

\begin{lem}\label{lem:basic}
Let $G$ and $H$ be groups acting on themselves and on each other. Then the following hold in $G\bt H$ for all $g\in G$ and $h\in H$:
\begin{align}
1\bt h &= 1 = g\bt 1, \tag{2.3.1}\\
{}^{g}g^{-1}\bt{}^{g}h &= (g\bt h)^{-1}. \tag{2.3.2}
\end{align}
In particular, every element of $G\bt H$ is a finite product of elements of $Y=\{g\bt h\mid g\in G,\ h\in H\}$.
\end{lem}

\begin{proof}
Setting $g=g'=1$ in (1.5.1) and noting that ${}^{1}1=1$ and ${}^{1}h=h$, we obtain $1\bt h=(1\bt h)(1\bt h)$, hence $1\bt h=1$. Similarly, setting $h=h'=1$ in (1.5.2) yields $g\bt 1=1$. To prove (2.3.2), expansion using (1.5.1) and (2.3.1) yields $1=gg^{-1}\bt h=({}^{g}g^{-1}\bt{}^{g}h)(g\bt h)$. The last assertion follows from (2.3.2).
\end{proof}

The following lemma is the key to the normality of $Y$. In the case of Inassaridze's tensor product, the normality of $Y$ is a consequence of the third relation (1.3.3) (see the proof of \cite[Theorem~6.3]{LT12}). Here we show that it already follows from (1.5.1) and (1.5.2), for arbitrary actions by automorphisms.

\begin{lem}\label{lem:identity}
Let $G$ and $H$ be groups acting on themselves and on each other. Then
\begin{equation}
(g\bt h)\Big({}^{({}^{h}g)}\big({}^{h}g'\big)\bt{}^{({}^{h}g)}\big({}^{h}h'\big)\Big)(g\bt h)^{-1}={}^{({}^{g}h)}\big({}^{g}g'\big)\bt{}^{({}^{g}h)}\big({}^{g}h'\big) \tag{2.4.1}
\end{equation}
for all $g,g'\in G$ and $h,h'\in H$.
\end{lem}

\begin{proof}
Let $g,g'\in G$ and $h,h'\in H$. Expansion using (1.5.1) and then (1.5.2), together with ${}^{g}(hh')={}^{g}h\,{}^{g}h'$, yields
\begin{equation}
gg'\bt hh'=\big({}^{g}g'\bt{}^{g}h\big)\Big({}^{({}^{g}h)}\big({}^{g}g'\big)\bt{}^{({}^{g}h)}\big({}^{g}h'\big)\Big)(g\bt h)\big({}^{h}g\bt{}^{h}h'\big). \tag{2.4.2}
\end{equation}
On the other hand, expansion using (1.5.2) and then (1.5.1), together with ${}^{h}(gg')={}^{h}g\,{}^{h}g'$, yields
\begin{equation}
gg'\bt hh'=\big({}^{g}g'\bt{}^{g}h\big)(g\bt h)\Big({}^{({}^{h}g)}\big({}^{h}g'\big)\bt{}^{({}^{h}g)}\big({}^{h}h'\big)\Big)\big({}^{h}g\bt{}^{h}h'\big). \tag{2.4.3}
\end{equation}
Comparing (2.4.2) and (2.4.3), we obtain (2.4.1).
\end{proof}

Using Lemma~\ref{lem:identity}, we now compute the conjugate of a tensor by a tensor and by its inverse. This shows that $Y$ is a normal set.

\begin{lem}\label{lem:normal}
Let $G$ and $H$ be groups acting on themselves and on each other. Let $a\in G$ and $b\in H$, and put $x={}^{b}a\in G$ and $y={}^{a}b\in H$. Then for all $u\in G$ and $v\in H$,
\begin{align}
(a\bt b)(u\bt v)(a\bt b)^{-1} &= {}^{y}\Big({}^{a}\big({}^{b^{-1}}({}^{x^{-1}}u)\big)\Big)\bt{}^{y}\Big({}^{a}\big({}^{b^{-1}}({}^{x^{-1}}v)\big)\Big), \tag{2.5.1}\\
(a\bt b)^{-1}(u\bt v)(a\bt b) &= {}^{x}\Big({}^{b}\big({}^{a^{-1}}({}^{y^{-1}}u)\big)\Big)\bt{}^{x}\Big({}^{b}\big({}^{a^{-1}}({}^{y^{-1}}v)\big)\Big). \tag{2.5.2}
\end{align}
In particular, $Y=\{g\bt h\mid g\in G,\ h\in H\}$ is a normal set. If $G$ and $H$ are finite, then $Y$ is a finite normal set.
\end{lem}

\begin{proof}
Throughout the proof, in each expression ${}^{c}w$ the action is determined by whether $c$ and $w$ lie in $G$ or in $H$. Since all the actions are homomorphisms into the automorphism groups, we have ${}^{c}\big({}^{c^{-1}}w\big)=w$ for all $c,w\in G\cup H$.

To prove (2.5.1), let $u\in G$ and $v\in H$, and apply (2.4.1) with $g=a$, $h=b$,
\[
g'={}^{b^{-1}}\big({}^{x^{-1}}u\big)\in G\qquad\text{and}\qquad h'={}^{b^{-1}}\big({}^{x^{-1}}v\big)\in H.
\]
Then ${}^{h}g={}^{b}a=x$ and ${}^{g}h={}^{a}b=y$. Moreover, ${}^{h}g'={}^{b}\big({}^{b^{-1}}({}^{x^{-1}}u)\big)={}^{x^{-1}}u$, and hence ${}^{({}^{h}g)}\big({}^{h}g'\big)={}^{x}\big({}^{x^{-1}}u\big)=u$. Similarly, ${}^{({}^{h}g)}\big({}^{h}h'\big)={}^{x}\big({}^{x^{-1}}v\big)=v$. Therefore the left hand side of (2.4.1) is $(a\bt b)(u\bt v)(a\bt b)^{-1}$, and the right hand side of (2.4.1) is
\[
{}^{y}\big({}^{a}g'\big)\bt{}^{y}\big({}^{a}h'\big)={}^{y}\Big({}^{a}\big({}^{b^{-1}}({}^{x^{-1}}u)\big)\Big)\bt{}^{y}\Big({}^{a}\big({}^{b^{-1}}({}^{x^{-1}}v)\big)\Big).
\]
This proves (2.5.1).

To prove (2.5.2), let $u\in G$ and $v\in H$, and put
\[
u_0={}^{x}\Big({}^{b}\big({}^{a^{-1}}({}^{y^{-1}}u)\big)\Big)\in G\qquad\text{and}\qquad v_0={}^{x}\Big({}^{b}\big({}^{a^{-1}}({}^{y^{-1}}v)\big)\Big)\in H.
\]
Then ${}^{x^{-1}}u_0={}^{b}\big({}^{a^{-1}}({}^{y^{-1}}u)\big)$, hence ${}^{b^{-1}}\big({}^{x^{-1}}u_0\big)={}^{a^{-1}}\big({}^{y^{-1}}u\big)$, hence ${}^{a}\big({}^{b^{-1}}({}^{x^{-1}}u_0)\big)={}^{y^{-1}}u$, and finally ${}^{y}\Big({}^{a}\big({}^{b^{-1}}({}^{x^{-1}}u_0)\big)\Big)=u$. Similarly, ${}^{y}\Big({}^{a}\big({}^{b^{-1}}({}^{x^{-1}}v_0)\big)\Big)=v$. Hence by (2.5.1) applied to $u_0\bt v_0$, we obtain $(a\bt b)(u_0\bt v_0)(a\bt b)^{-1}=u\bt v$. Conjugating both sides by $(a\bt b)^{-1}$ gives $(a\bt b)^{-1}(u\bt v)(a\bt b)=u_0\bt v_0$, which is (2.5.2).

By (2.5.1) and (2.5.2), conjugation by $a\bt b$ and by $(a\bt b)^{-1}$ maps $Y$ into $Y$, for all $a\in G$ and $b\in H$. By Lemma~\ref{lem:basic}, every element of $G\bt H$ is a finite product of elements of $Y$. By induction on the number of factors, it follows that $z\,Y\,z^{-1}\subseteq Y$ for all $z\in G\bt H$. Hence $Y$ is a normal set. If $G$ and $H$ are finite, then obviously $Y$ is finite.
\end{proof}

\section{Finiteness of $G\bt H$ and of Inassaridze's tensor product}

In this section, we prove that $G\bt H$ is finite whenever $G$ and $H$ are finite, and we deduce the finiteness of Inassaridze's tensor product. We start with a lemma, which is the heart of the proof. It concerns maps into abelian groups satisfying conditions of a crossed pairing in the sense of Remark 3 of \cite{BJR87}. The lemma requires no compatibility condition; in fact only the bijectivity of the actions is used.

\begin{lem}\label{lem:avg}
Let $G$ and $H$ be groups acting on themselves and on each other. Let $A$ be an abelian group, written additively, and let $E\colon G\times H\to A$ be a map such that
\begin{align}
E(gg',h) &= E\big({}^{g}g',{}^{g}h\big)+E(g,h), \tag{3.1.1}\\
E(g,hh') &= E(g,h)+E\big({}^{h}g,{}^{h}h'\big), \tag{3.1.2}
\end{align}
for all $g,g'\in G$ and $h,h'\in H$. If $G$ and $H$ are finite, then $|G||H|\,E(g,h)=0$ for all $g\in G$ and $h\in H$.
\end{lem}

\begin{proof}
For $h\in H$ and $g\in G$ set
\[
L(h)=\sum_{a\in G}E(a,h)\qquad\text{and}\qquad R(g)=\sum_{b\in H}E(g,b).
\]
Fix $g\in G$ and $h\in H$. Since $\rho_G(g)\in\Aut(G)$, as $g'$ runs over $G$ so do $gg'$ and ${}^{g}g'$. Hence summing (3.1.1) over all $g'\in G$ yields $L(h)=L({}^{g}h)+|G|\,E(g,h)$, that is
\begin{equation}
|G|\,E(g,h)=L(h)-L\big({}^{g}h\big). \tag{3.1.3}
\end{equation}
Similarly, summing (3.1.2) over all $h'\in H$ yields,
\begin{equation}
|H|\,E(g,h)=R(g)-R\big({}^{h}g\big). \tag{3.1.4}
\end{equation}
Now fix $h\in H$ and sum (3.1.4) over all $g\in G$. Since $\sigma_H(h)\in\Aut(G)$, the map $g\mapsto{}^{h}g$ is a bijection of $G$, and therefore $\sum_{g\in G}R({}^{h}g)=\sum_{g\in G}R(g)$. This leads to
\begin{equation}
|H|\,L(h)=0\qquad\text{for all } h\in H. \tag{3.1.5}
\end{equation}
Finally, multiplying (3.1.3) by $|H|$ and using (3.1.5) for $h$ and for ${}^{g}h$, we obtain
\[
|G||H|\,E(g,h)=|H|\,L(h)-|H|\,L\big({}^{g}h\big)=0,
\]
the desired result.
\end{proof}

Next we will prove that the abelianization $(G\bt H)^{ab}$ is finite, provided $G$ and $H$ are finite.

\begin{prop}\label{prop:ab}
Let $G$ and $H$ be finite groups acting on themselves and on each other. Then
\[
(g\bt h)^{|G||H|}\in[G\bt H,G\bt H]
\]
for all $g\in G$ and $h\in H$. In particular, the abelianization $(G\bt H)^{ab}$ is a finite abelian group whose exponent divides $|G||H|$.
\end{prop}

\begin{proof}
Let $\pi\colon G\bt H\to(G\bt H)^{ab}$ be the canonical epimorphism, and write $(G\bt H)^{ab}$ additively. Define $E\colon G\times H\to(G\bt H)^{ab}$ by $E(g,h)=\pi(g\bt h)$. Applying $\pi$ to (1.5.1) and (1.5.2) shows that $E$ satisfies (3.1.1) and (3.1.2). By Lemma~\ref{lem:avg}, we obtain $\pi\big((g\bt h)^{|G||H|}\big)=|G||H|\,E(g,h)=0$, that is $(g\bt h)^{|G||H|}\in\ker\pi=[G\bt H,G\bt H]$. Since $(G\bt H)^{ab}$ is generated by the finitely many elements $E(g,h)$, $g\in G$, $h\in H$, each of order dividing $|G||H|$, the last part of our claim follows.
\end{proof}

Now we are ready to prove the main theorem.

\begin{thm}\label{thm:main}
Let $G$ and $H$ be finite groups acting on themselves and on each other. Then $G\bt H$ is finite.
\end{thm}

\begin{proof}
By Lemma~\ref{lem:normal}, $Y=\{g\bt h\mid g\in G,\ h\in H\}$ is a finite normal set which generates $G\bt H$. By Proposition~\ref{prop:ab}, the abelianization $(G\bt H)^{ab}$ is finite. Hence by Lemma~\ref{lem:B}, $G\bt H$ is finite.
\end{proof}

As a consequence, we obtain the finiteness of Inassaridze's tensor product.

\begin{cor}\label{cor:inas}
Let $G$ and $H$ be finite groups acting on themselves by conjugation and each of which acts on the other, not necessarily compatibly. Then $G\otimes H$ is finite.
\end{cor}

\begin{proof}
Let $G\bt H$ be the group of Definition~\ref{def:box}, where $G$ and $H$ act on themselves by conjugation. Then (1.5.1) and (1.5.2) are precisely (1.3.1) and (1.3.2). Hence the assignment $g\bt h\mapsto g\otimes h$ defines an epimorphism $G\bt H\to G\otimes H$, whose kernel is the normal subgroup generated by the relations (1.3.3). By Theorem~\ref{thm:main}, $G\bt H$ is finite, and hence so is $G\otimes H$.
\end{proof}

\begin{cor}\label{cor:box}
Let $G$ and $H$ be finite groups acting on each other fully compatibly. Then the box-tensor product $G\boxt H$ is finite. Let $G$ and $H$ be finite groups acting on themselves by conjugation and on each other compatibly. Then the nonabelian tensor product $G\otimes H$ is finite.
\end{cor}

\begin{proof}
The first statement is the special case of Theorem~\ref{thm:main} for fully compatible actions. The second statement is the special case of Corollary~\ref{cor:inas} for compatible actions.
\end{proof}

\begin{rem}\label{rem:main}
(a) The first part of Corollary~\ref{cor:box} is \cite[Theorem~4.8]{LT12}, and the second part is Ellis's theorem \cite{Ell87}. Theorem~\ref{thm:main} also shows that the groups $G\bt H$ considered in \cite[Example~2.2]{LT12} for the actions defined by $\psi_b$ and $\psi_a$, which are not fully compatible, are finite.

(b) The third relation (1.3.3) of Inassaridze's tensor product plays no role in the proof of Corollary~\ref{cor:inas}.
\end{rem}

The next example shows that Corollary~\ref{cor:inas} applies in situations which are not covered by \cite[Theorem~6.4]{LT12}.

\begin{example}\label{ex:V4}
Let $G\cong H\cong V_4=\{a,b,ab,e\}$, the Klein four group, acting on themselves by conjugation, that is trivially. Let $\alpha\in\Aut(V_4)$ be the automorphism which interchanges $a$ and $b$ and fixes $ab$. Define $\sigma_G\colon G\to\Aut(H)$ and $\sigma_H\colon H\to\Aut(G)$ by
\[
\sigma_G(a)=\sigma_G(ab)=\alpha,\quad\sigma_G(b)=1,\qquad\sigma_H(a)=\sigma_H(ab)=\alpha,\quad\sigma_H(b)=1.
\]
Since $G$ and $H$ are abelian and act trivially on themselves, the element $hgh^{-1}\in G\ast H$ acts on $H$ as $g$ does. Hence (1.1.1) reads ${}^{({}^{h}g)}h'={}^{g}h'$. Taking $g=h=h'=a$, we obtain ${}^{({}^{a}a)}a={}^{b}a=a$, whereas ${}^{a}a=b$. Thus (1.1.1) does not hold, and by symmetry neither does (1.1.2). So the mutual actions are not half compatible, and \cite[Theorem~6.4]{LT12} does not apply, but $G\otimes H$ is finite, by Corollary~\ref{cor:inas}.
\end{example}

\section{Applications to low dimensional nonabelian homology groups}

In this section we deduce the finiteness of the nonabelian homology groups $H_0(G,A)$ and $H_1(G,A)$ introduced in \cite{Ina96a}. Let $G$ and $A$ be groups acting on themselves by conjugation and each of which acts on the other. Recall from \cite{Ina96a} that
\[
H_1(G,A)=\ker\lambda'\qquad\text{and}\qquad H_0(G,A)=\coker\lambda',
\]
where $\lambda'\colon G\otimes A\to A/A'$ is given by $\lambda'(g\otimes a)={}^{g}aa^{-1}A'$, and $A'$ is the normal subgroup of $A$ generated by the elements ${}^{({}^{a}g)}a'\,{}^{aga^{-1}}a'^{-1}$ for all $a,a'\in A$ and $g\in G$. Here $G\otimes A$ denotes Inassaridze's tensor product. Note that $A'$ is trivial when the mutual actions are compatible.

\begin{cor}\label{cor:hom}
Let $G$ and $A$ be finite groups acting on each other. Then the nonabelian homology groups $H_i(G,A)$ are finite when $i=0,1$.
\end{cor}

\begin{proof}
By Corollary~\ref{cor:inas}, $G\otimes A$ is finite, and hence so is $H_1(G,A)=\ker\lambda'$. Since $A$ is finite, so is $H_0(G,A)=\coker\lambda'$.
\end{proof}


\begin{thebibliography}{99}
\bibitem{BJR87} R. Brown, D.L. Johnson, E.F. Robertson, Some computations of non-abelian tensor products of groups, J. Algebra 111 (1987) 177--202.
\bibitem{BL84} R. Brown, J.-L. Loday, Excision homotopique en basse dimension, C. R. Acad. Sci. Paris S\'er. I Math. 298 (1984) 353--356.
\bibitem{BL87} R. Brown, J.-L. Loday, Van Kampen theorems for diagrams of spaces, Topology 26 (1987) 311--335.
\bibitem{BNR21} R. Bastos, I.N. Nakaoka, N.R. Rocco, Finiteness conditions for the box-tensor product of groups and related constructions, J. Algebra 587 (2021) 594--612.
\bibitem{BNR25} R. Bastos, I.N. Nakaoka, N.R. Rocco, On the finiteness of the non-abelian tensor product of groups, J. Algebra 664 (2025) 251--267.
\bibitem{Die37} A.P. Dietzmann, \"Uber $p$-Gruppen, Dokl. Akad. Nauk SSSR 15 (1937) 71--76.
\bibitem{Ell87} G.J. Ellis, The non-abelian tensor product of finite groups is finite, J. Algebra 111 (1987) 203--205.
\bibitem{Ina96a} N. Inassaridze, Non-abelian tensor products and non-abelian homology of groups, J. Pure Appl. Algebra 112 (1996) 191--205.
\bibitem{Ina96b} N. Inassaridze, Finiteness of a nonabelian tensor product of groups, Theory Appl. Categ. 2 (1996) 55--61.
\bibitem{Kap99} L.-C. Kappe, Nonabelian tensor products of groups: the commutator connection, in: Groups St. Andrews 1997 in Bath, II, London Math. Soc. Lecture Note Ser., vol. 261, Cambridge Univ. Press, Cambridge, 1999, pp. 447--454.
\bibitem{LT12} M. Ladra, V.Z. Thomas, Two generalizations of the nonabelian tensor product, J. Algebra 369 (2012) 96--113.
\bibitem{Rob82} D.J.S. Robinson, A Course in the Theory of Groups, Grad. Texts in Math., vol. 80, Springer-Verlag, New York, Berlin, 1982.
\bibitem{Tho10} V.Z. Thomas, The non-abelian tensor product of finite groups is finite: A homology-free proof, Glasg. Math. J. 52 (2010) 473--477.
\end{thebibliography}
\end{document}